\documentclass[reqno]{amsart}
\usepackage{graphicx} 
\usepackage{amssymb, amsmath, amsthm,xcolor}
\usepackage{tikz}
\usepackage{pgfplots}
\usepackage{url}
\usepackage{subcaption}

\theoremstyle{plain}
\newtheorem{thm}{Theorem}

\theoremstyle{definition}

\newtheorem{ex}[thm]{Example}

\newtheorem{re}[thm]{Remark}

\DeclareMathOperator{\RR}{\mathbb{R}}

\newcommand{\ve}[1]{\mathbf{#1}}%

\title{Global Positioning on Earth}

\author{Mireille Boutin} \address{Technical University of Eindhoven,
  5600 MB Eindhoven, Netherlands} \email{m.boutin@tue.nl}%
\author{Rob Eggermont} \address{Technical University of Eindhoven,
  5600 MB Eindhoven, Netherlands} \email{r.h.eggermont@tue.nl}
\author{Gregor Kemper} \address{Technical University of Munich,
  Germany; TUM School of Computation, Information and Technology,
  Department of Mathematics} \email{kemper@tum.de}


\begin{document}

\maketitle
\begin{abstract}
  Contrary to popular belief, the global positioning problem on earth
  may have more than one solutions even if the user position is
  restricted to a sphere.  With 3 satellites, we show that there can
  be up to 4 solutions on a sphere.  With 4 or more satellites, we
  show that, for any pair of points on a sphere, there is a family of
  hyperboloids of revolution such that if the satellites are placed on
  one sheet of one of these hyperboloid, then the global positioning
  problem has both points as solutions.
  We give solution methods that yield the correct number of solutions
  on/near a sphere.
\end{abstract}

\section{Introduction }
It is a long standing assumption that the global positioning (GP) problem has a unique solution in the case where the user to be located is on the surface of the earth and the satellites are in orbit around it. For example, Abel and Chaffee \cite{Abel:Chaffee:1991}, one of seminal references on GP, states that: ``A unique solution is vitrually (sic) guaranteed for users operating near the earth.'' This is made precise in Figure 5 of their paper, in which the earth is pictured to be entirely within the region in which a unique solution exists. On the website of the European Space Agency~\cite{ESA} we find: ``The other solution is far from the earth surface.''

As we show in Sect.~\ref{sec:two_solutionsa_on_earth}, this assumption is incorrect. However, the fact that the user position is close to the earth allows us to formulate different solution methods (Sect.~\ref{sec:solution_onsphere}). In particular, a sphere assumption allows for location with only $3$ satellites, but in that case there can be up to $4$ solutions. 
For 4 or more satellites, we propose a simple numerical procedure in which two initial guesses on a sphere are refined numerically so to find a user position assumed to be near a sphere.   Numerical experiments show how our proposed method allows us to correctly obtain either $1$ or $2$ solutions near a sphere. 
For simplicity, our results are phrased in dimension $n=3$ but do carry over to any dimension $n\geq 2$.

\subsection{Problem Statement}
We define the global positioning (GP) problem as follows.
A user at an unknown position $\ve x \in \RR^3$ receives a signal from $m$ satellites at known positions $\ve s_1,\ldots, \ve s_m\in \RR^3$. The signal from satellite $i$ contains its position $\ve s_i$ at the time of emission $\tau_i$, along with $\tau_i$. The user receives the signal $(\ve s_i, \tau_i)$ at a time $\bar{\tau}_i$, as measured by their own clock.
The satellite clocks are precisely synchronized together but not with the user clock. 
Let $t_i$ be the difference between the emission time $\tau_i$ and the (local) reception time
$\bar{\tau}_i$: $t_i=\bar{\tau}_i -\tau_i$, for $i=1,\ldots, m$. Then there exists an unknown offset $t\in \RR$ such that $t_i -t$ is the true time taken by the signal to travel from satellite $i$ to the user. 
The value of $t_i -t$ is proportional to the distance between the satellite position $\ve s_i$ and the user position $\ve x$.
Using units of time so that the signal transmission speed is $1$, we thus have
\begin{equation}
\|\ve s_i -\ve x \| = t_i -t \geq 0, \quad i=1,\ldots, m. 
\label{eq:GPS_equations}
\end{equation}
As shown in \cite{boutin2024global}, this system of equations can be rewritten as \begin{equation} 
\label{eqLin}%
      \begin{pmatrix}
    -2 t_1 & 2 \ve s_1^T & -1 \\
    \vdots & \vdots & \vdots \\
    -2 t_m & 2 \ve s_m^T & -1
  \end{pmatrix} 
      \begin{pmatrix}
        t \\ \ve x \\ \| \ve x\|^2- t^2
      \end{pmatrix} =
      \begin{pmatrix}
        \lVert\ve s_1\rVert^2 - t_1^2 \\ \vdots \\ \lVert\ve
        s_m\rVert^2 - t_m^2
      \end{pmatrix}\!\!,\ \text{and\ } t_i>t, \text{ }i=1,\ldots,
      m. 
    \end{equation}
    Assuming that the user is located on a sphere of radius $r$ centered at  $\ve c$ adds the constraint $\| \ve x - \ve c\|= r $.

\section{Two solutions on a sphere}
\label{sec:two_solutionsa_on_earth}

Both the earth and the orbits of satellites are nearly circular \cite{hofmann2012global}. As previously mentioned, this fact has been wrongly assumed to make the GP solution unique. 
This section goes into more details to explain how things can go wrong.
According to Thm.~3.1 of \cite{boutin2024global}, for any $m\geq 4$, the general GP problem (not on a sphere) has two distinct solutions $\ve x$ and $\ve x'$ if and only if the satellites lie on one sheet of a hyperboloid of revolution with focal points $\ve x, \ve x'$; in all other cases, the problem has only one solution. Conversely, given $\ve x, \ve x'$ we can classify all satellite configurations for which the general GP problem can have solutions $\ve x, \ve x'$ by finding all hyperboloids of revolutions that have these points as focal points. The following theorem describes these hyperboloids.

\begin{thm} Let $\ve x \neq \ve x' \in \RR^3$. Let $\ve m = \frac{1}{2}(\ve x + \ve x')$, let $c = \| \ve x - \ve x'\|$ and let $\tilde{\ve u} = \frac{1}{c}(\ve x - \ve x')$. Then for every $0 < a < c$, the hyperboloid of revolution given by the following equation  has focal points $\ve x, \ve x'$: \begin{equation}
(\ve p - \ve m)^T((\frac{1}{a^2} + \frac{1}{c^2-a^2})\tilde{\ve u} \tilde{\ve u}^T - \frac{1}{c^2-a^2}\mathrm{Id}_3)(\ve p - \ve m) = 1.
\label{eq:family}
\end{equation}
\end{thm}

\begin{figure}
\begin{subfigure}{.45\textwidth}
 \centering
  \includegraphics[width=.8\linewidth]{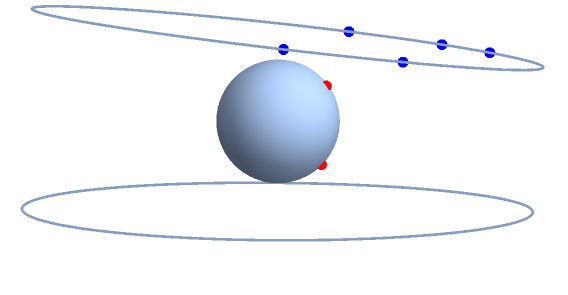}
\end{subfigure}\hfill
\begin{subfigure}{.45\textwidth}
  \centering
  \includegraphics[width=.4\linewidth]{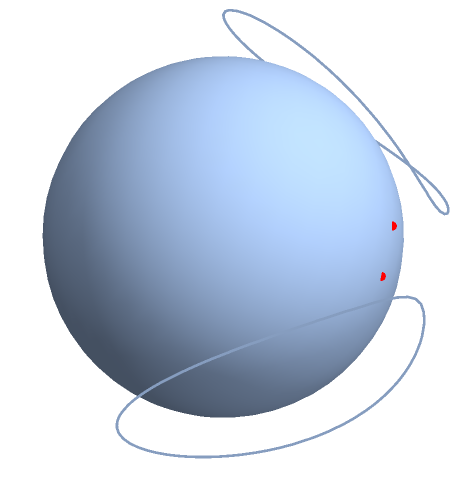}
\end{subfigure}
\caption{(Left) Two solutions (red dots) to the GP problem on a sphere and one possible locus of forbidden satellite positions. (Right) Visibly non-circular loci.}
\label{fig:2solutions_sphere}
\end{figure}

A hyperboloid of revolution with Eq.~\ref{eq:family} has axis of symmetry $t\tilde{\ve u} + \ve m$ with $t$ a parameter. For one sheet of such a hyperboloid, the points $\ve x, \ve x'$ are indistinguishable in the sense that the distance from $\ve x$ to the sheet equals that from $\ve x'$ to the sheet plus a fixed constant. Thus there is always 
two possible offsets, regardless of the number of satellites on this sheet. Given such a hyperboloid of revolution, we can find problematic satellite positions by intersecting the hyperboloid with the satellite orbit
and look at the connected components of this intersection. Fig.~\ref{fig:2solutions_sphere} shows an example of a bad configuration.

More generally, one could ask how often it happens that all visible satellites lie on such a sheet, i.e.~allow for bad configurations with multiple solutions. Eq.~\ref{eq:family} has $7$ free parameters ($\ve x, \ve x'$, $a$), so one should expect that when at most $7$ satellites are in view, there are some quadrics of revolution passing through them, and we expect infinitely many (including hyperboloids of revolution) with fewer satellites. The article \cite{gferrer2009quadrics} checks this in the case where the points are not cospherical. To summarize, when there are $7$ or fewer satellites in general position, there typically are multiple quadrics of revolution passing through them, with more options when there are fewer satellites. While these results do not apply in our situation since satellites orbit around the earth,
it is at least plausible that positions as described in the above theorem may actually be found in practice, especially when the user views few satellites.

Many GP systems use iterative solution methods that start from an initial guess. Thus, they return a single solution, even if there are two. Given one (exact) solution of the system, one might want to check whether the satellite and user positions are in a bad configuration (i.e., with two solutions). This can be done by using the arrival times and satellite positions of at least four satellite signals to compute the vector $\ve u$ as described in Section \ref{sec:near_sphere}.  If $|\ve u|\leq 1$, then this system of equations only has one solution. 
If $|\ve u|>1$, then there might be two, and the second one could potentially be on the sphere (earth) as well. There are two ways to determine if there is another solution on the sphere: either compute all the solutions using the method of \cite{boutin2024global} and check if they are on the sphere (or on earth), or compute all the solutions on the sphere, as described in Sect.~\ref{sec:solution_onsphere}.

\section{Global positioning on or near a sphere}
\label{sec:solution_onsphere}

\subsection{On a Sphere With Three satellites}
\label{sec:three_satellites}
\newcommand{\upto}{,\ldots,}%

For getting a GP fix, a minimum of four satellites in view is
required. However, using the information that the user is located on
the surface of the earth, one might hope that this minimum number is
reduced to three, or at least three satellites in view may lead to
finitely many solutions. The following result fulfills this hope if the satellite positions are not collinear.
It shows that there are at most~$4$ solutions for~$t$. Having a solution
for~$t$, this can be substituted into the original problem, which then becomes
a problem of intersecting four spheres in $\RR^3$. But this problem is well known to have
at most two solutions.

\begin{thm} \label{ThreeSatellites}
  Let $\ve c,\ve s_1,\ve s_2,\ve s_3 \in \RR^3$ be points such that $\ve
  s_i$ are not collinear, and let
  $r,t_1,t_2,t_3 \in \RR$ be numbers with $r \ge
  0$. Consider the set
  \[
    {\mathcal S} := \bigl\{(\ve x,t) \in \RR^3 \times \RR \mid \lVert\ve x - \ve
    c\rVert = r \ \text{and} \ \lVert\ve x - \ve s_i\rVert = |t_i
    - t|\ \text{for} \ i = 1 \upto n\bigr\}
  \]
  of solutions of the GP problem (slightly loosened by using the
  absolute value of $t_i - t$), with~$\ve x$ restricted to the sphere
  given by the first equation.  Then there is a polynomial~$f$ of
  degree four such that every $(\ve x,t) \in \mathcal S$ satisfies
  $f(t) = 0$. One can calculate~$f$ by setting $f := \det(C)$ with $C$
  given by~\eqref{eqCM} below.
\end{thm}

\begin{proof}
  Let $(\ve x,t) \in \mathcal S$. The Cayley-Menger matrix $C$ of any
  collection of points in $\RR^3$ has rank at most~$5$, and in
  particular this applies to the points $\ve x$, $\ve c$,
  $\ve s_1,\ve s_2,\ve s_3$. Since $(\ve x,t) \in \mathcal S$, the
  Cayley-Menger matrix for these points is, with
  $d_i := \lVert\ve s_i - \ve c\rVert$ and
  $d_{i,j} := \lVert\ve s_i - \ve s_j\rVert$,
  \begin{equation} \label{eqCM}%
    C =
    \begin{pmatrix}
      0 & 1 & 1 & 1 & 1 & 1 \\
      1 & 0 & r^2 & (t_1 -t)^2 & (t_2 - t)^2 & (t_3
      - t)^2 \\
      1 & r^2 & 0 & d_1^2 & d_2^2 & d_3^2 \\
      1 & (t_1 - t)^2 & d_1^2 & 0 & d_{1,2}^2 & d_{1,3}^2 \\
      1 & (t_2 - t)^2 & d_2^2 & d_{2,1}^2 & 0 & d_{2,3}^2 \\
      1 & (t_3 - t)^2 & d_3^2 & d_{3,1}^2 & d_{3,2}^2 & 0
    \end{pmatrix} \in \RR^{6 \times 6}.
  \end{equation}
  So $\det(C) = 0$. Consider $f(t) := \det(C)$ as a function of~$t$.
  As such, it is clearly a polynomial 
  in~$t$. To show that~$f(t)$ has degree~$4$,
  we form the rational function $g(t) := t^4 f(1/t)$ . If we can show that~$g$ is 
  actually a polynomial and $g(0) \ne 0$, then indeed~$f(t)$ has degree~$4$. To 
  compute~$g(t)$ we perform the multiplication with $t^4$ be multiplying the second
  row and the second column of $C$, with~$t$ substituted by~$1/t$, by~$t^2$. It is then easy to see that $g(t)$ is a polynomial in~$t$, and that
  \[
  g(0) = 
    -\det
  \begin{pmatrix}
      0 & 1 & 1 & 1 \\
      1 & 0 & d_{1,2}^2 & d_{1,3}^2 \\
      1 & d_{2,1}^2 & 0 & d_{2,3}^2 \\
      1 & d_{3,1}^2 & d_{3,2}^2 & 0
    \end{pmatrix},
  \]
  which
  is the negative of the Cayley-Menger determinant of the points~$\ve s_1,\ve s_2,\ve s_3$
  and thus nonzero since the~$\ve s_i$ are assumed not collinear.
\end{proof}

\begin{re}
        Theorem~\ref{ThreeSatellites} shows that there are at most four solutions
        for~$t$, and by what we have said before the theorem, each of these can lead to at most two
        solutions of the GP problem, giving a total of $|\mathcal S| \le 8$ solutions.
        In fact, it can be shown that there cannot be more than four
        solutions in total. For once, if the four points $\ve c,\ve s_1,\ve s_2,\ve s_3$ are not coplanar,
        then for each zero~$t$ of~$f$ there exists exactly one~$\ve x \in \RR^3$ such that $(\ve x,t) \in \mathcal S$. 
         If, on the other hand, the four points are coplanar, there are 
        at most two points~$\ve x$ with $(\ve x,t) \in \mathcal S$, but there may also be none. But it can be shown that also in the coplanar case, the
        total number $|\mathcal S|$ of solutions does not exceed four.
        
        We also have an argument that takes care
        of the case in which the earth is modeled as a spheroid rather than a sphere. 
        The results are the same.
\end{re}

\begin{ex}
    It is quite rare to find a configuration where there are really four solutions that solve the original GP problem~\eqref{eq:GPS_equations}, rather than only the squared equations. We give two examples. The computations were done using Magma~\cite{magma}.
    \begin{enumerate}
    \renewcommand{\labelenumi}{(\theenumi)}
    \renewcommand{\theenumi}{\arabic{enumi}}
        \item A 3D example is given by $\ve s_1 = (-4,6,6)$, $\ve s_2 = (0,1,2)$, $\ve s_3 = (-1,5,9)$, $t_1 = -2$, $t_2 = -9$, $t_3 = -1$.
        From this we get the polynomial $f = 4232 t^4 + 188416 t^3 + 
        3141776 t^2 + 23254272 t + 64463912$, and its approximate roots are $-11.8922$, $-11.7298$,$-10.4779$, and $-10.4216$. It is easy to confirm that for each of these roots~$t$ there exists exactly one solution $(\ve x,t)$ of the system~\eqref{eq:GPS_equations}.
    \item The next example is in two dimensions, and we present it graphically in Figure~\ref{fig:Two_sats} instead of giving precise coordinates. Being in dimension two means that there are two satellites, and we solve~\eqref{eq:GPS_equations} for points on a given circle. In Figure~\ref{fig:Two_sats}, the green circles around the satellites~$\ve s_i$ are tangent to the black circle and the red circles, which implies: (radius of green circle) + (radius of black or red circle)  = (distance between the corresponding satellite and solution point).
    So if the~$t_i$ are set to be the radii of the green circles, this shows that indeed the~$\ve x_i$ solve the GP problem.

    \begin{figure}[h]
        \centering
        \includegraphics[width=0.75\linewidth]{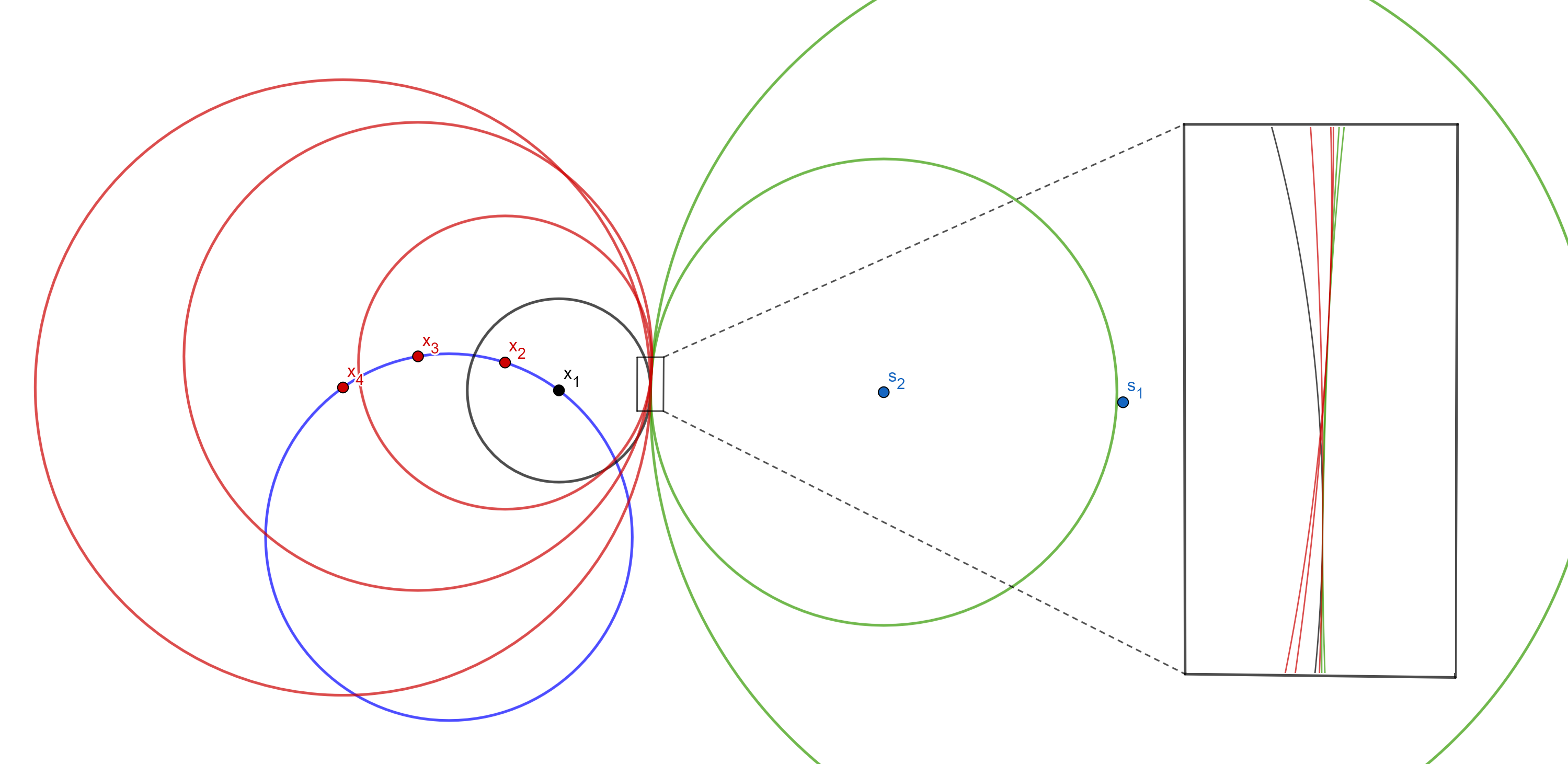}
        \caption{The solutions~$\ve x_i$ are on the blue circle, $\ve x_1$ being the ``true'' one. The enlargement shows that the green circles around the satellites~$\ve s_i$ are tangent to every one of the black and red circles around the~$\ve x_i$.}
        \label{fig:Two_sats}
    \end{figure}
    \end{enumerate}
\end{ex}

\subsection{On a Sphere With Four or More Satellites}
\label{sec:on_sphere_with4}
To solve the GP problem on a sphere with $m\geq 4$  satellites, we write eq.~(\ref{eqLin}) as 
\begin{eqnarray*} 
     B 
      \begin{pmatrix}
       \ve x \\ \| \ve x\|^2- t^2
      \end{pmatrix} =
      t  \begin{pmatrix}
    2 t_1 \\
    \vdots  \\
    2 t_m
  \end{pmatrix} +
      \begin{pmatrix}
        \lVert\ve s_1\rVert^2 - t_1^2 \\ \vdots \\ \lVert\ve
        s_m\rVert^2 - t_m^2
      \end{pmatrix}, \text{ where } B=\begin{pmatrix}
    2 \ve s_1^T & -1 \\
  \vdots & \vdots \\
    2 \ve s_m^T & -1
  \end{pmatrix} .
    \end{eqnarray*}
Assuming that four of the satellites are non-coplanar, then $B$  has full rank. Thus its Moore-Penrose pseudo inverse is $B^+=(B^T B)^{-1} B^T$ and we have  
\[
      \begin{pmatrix}
       \ve x \\ \| \ve x\|^2- t^2
      \end{pmatrix} =
      t   \begin{pmatrix}
    \ve u \\
    2 \alpha
  \end{pmatrix} + \begin{pmatrix}
    \ve v \\
 \beta
\end{pmatrix},
\]
where
\[
  \begin{pmatrix}
    \ve u \\
    2 \alpha
  \end{pmatrix}= B^+ \begin{pmatrix}
    2 t_1 \\
    \vdots  \\
    2 t_4
  \end{pmatrix},  
  \begin{pmatrix}
    \ve v \\
 \beta
  \end{pmatrix}= B^+  \begin{pmatrix}
        \lVert\ve s_1\rVert^2 - t_1^2 \\ \vdots \\ \lVert\ve
        s_4\rVert^2 - t_4^2
      \end{pmatrix}.
\]
Focusing on the equation for $\ve x$, we have $\ve x=\ve ut+\ve v$ which can be replaced in the constraint $\| \ve x - \ve c\|^2 -r^2=0 $ to yields the quadratic equation in $t$:
\begin{equation}
    \|\ve u \|^2 t^2 + 2 \ve u \cdot (\ve v - \ve c )+ \| \ve v - c\|^2-r^2=0,
\label{eq:quadratic_sphere}
\end{equation} 
which has $\leq 2$ solutions (guaranteed to be on the circle.)
Incorrect solutions can be eliminated by checking if they satisfy Eq.~\ref{eq:GPS_equations}.

\subsection{Near a Sphere}
\label{sec:near_sphere}
When the user is known to be located near not necessarily on, a sphere with known radius and center, then it is not possible to determine its location with $3$ satellites. With $4$ or more satellites, and without putting any constraint on the location of the user, there are at most $2$ solutions, as described in \cite{boutin2024global}. There are many solution methods in such scenarios. A straightforward approach is to compute the total least squares solution of Eq.~\ref{eqLin}, treating the vector of unknowns as if it had unrelated components. This approach only produces one solution. However, it is well suited to when the observed times of arrival $t_i$ are noisy. Moreover, the solution obtained can be refined using an iterative method. 
In practice, this often works quite well. However, the results are unpredictable in the vicinity of satellite configurations for which two solutions exist, as we show in our numerical experiments below.

Eq.~\ref{eq:quadratic_sphere} suggests an alternative strategy when the user is known to be located in the vicinity of the surface of a sphere with known center and radius: use the $2$ solutions of Eq.~\eqref{eq:quadratic_sphere} as initial guesses for an iterative method. This has the advantage of allowing for $2$ solutions to be found; in case where the solution is unique, one would expect both initial guesses to converge to the same solution.  
Thus this strategy is better suited for dealing with situations where the GP problem either has $2$ solutions or is near one with  $2$ solutions. 

\subsubsection{Numerical experiments}
\label{sec:experiments}
We compare $3$ solution approaches to find the position of a user on/near a sphere of known center and radius. 
The first one, {\em Iterative Least Squares} (ILS), uses the least squares approach of \cite{ESA_ILS} to refine (20 steps) the initial guess obtained from the total least squares solution of Eq.~(\ref{eqLin}). The second one,  {\em Solution on Sphere} (SoS), consists in computing the $2$ solutions of Eq.~\ref{eq:quadratic_sphere}. 
For the third method, {\em Refined Solution on Sphere} (RSoS),  we used those $2$ solutions as initial guesses for the iteration procedure of \cite{ESA_ILS} (20 steps). 

\begin{figure}[h]
\begin{subfigure}{.48\textwidth}
  \includegraphics[width=\linewidth]{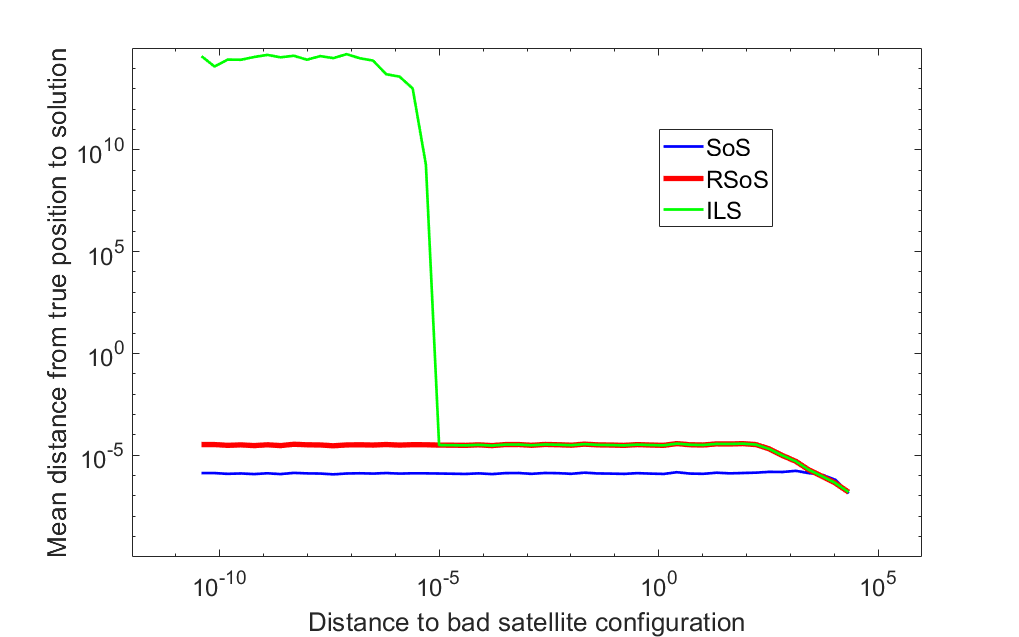}
\end{subfigure}\hfill
\begin{subfigure}{.48\textwidth}
  \includegraphics[width=\linewidth]{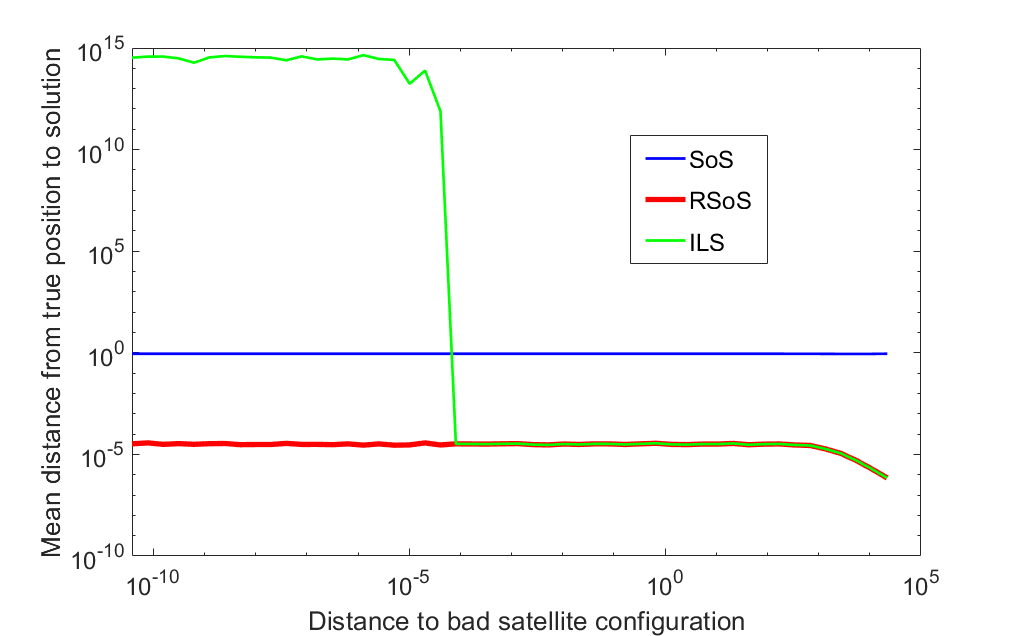}
\end{subfigure}
\caption{Mean distance to user position for ILS, SoS and RSoS 
when user is on the earth sphere (left) and $\sim 6$km above the earth (right).}
\label{fig:mean_distance_to_true_solution}
\end{figure}

\begin{figure}[h]
\begin{subfigure}{.48\textwidth}
  \includegraphics[width=\linewidth]{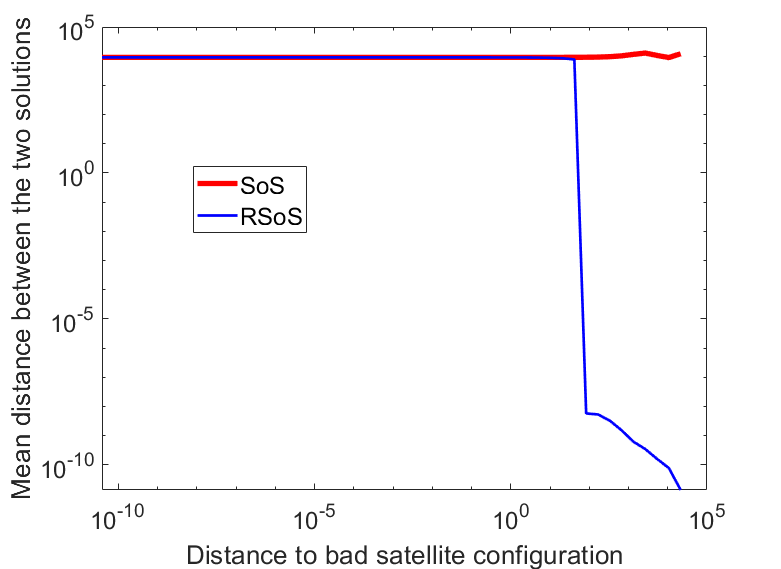}
\end{subfigure}\hfill
\begin{subfigure}{.48\textwidth}
  \includegraphics[width=\linewidth]{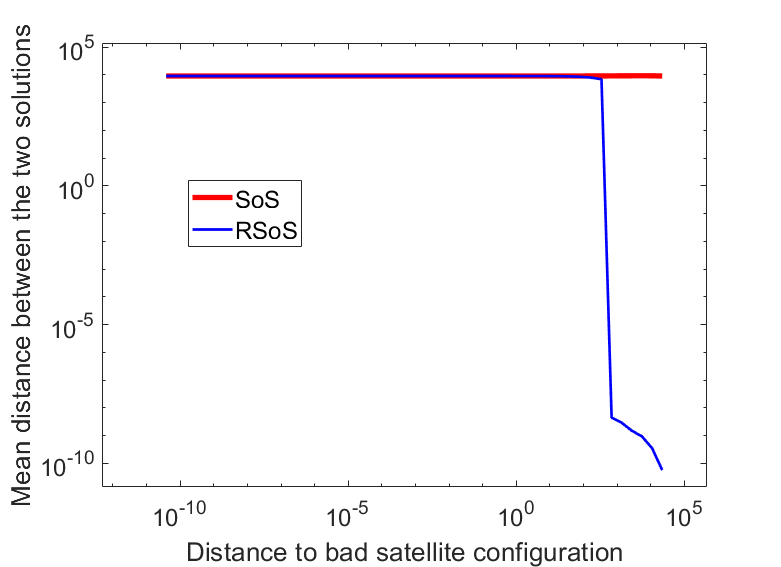}
\end{subfigure}
\caption{Mean distance between the two solutions of SoS and RSoS 
when the user is on the earth sphere (left) and $\sim 6$km above the earth sphere (right).}
\label{fig:distance_between_slns}
\end{figure}

In the first experiment, we put a user on a sphere with a 6400 km radius (approx.~the earth.) 
We put 5 satellites in a ``bad" configuration within a 26400 km circular orbit  of the sphere (approx.~GPS orbit.) 
In other words, the satellites were placed so to create a GP problem with two solutions on the user sphere.
We subsequently picked a random configuration of 5 satellites and constructed a path from the random configuration to the bad configuration. At each step along the path, the times of arrival for the signals of the 5 satellites were perturbed with additive zero mean Gaussian noise ($\sigma=10^{-8}$); 200 trials were performed at each step and the mean distance between the solution obtained and the correct one was computed. With ILS, only one solution was obtained each time; the average of the results are shown in Fig.~\ref{fig:mean_distance_to_true_solution} (left).
As one can see, the accuracy of the method is good (about $10^{-5}$ km) until the satellite configuration gets too close to the bad one, at which point the solution diverges.  With SoS and RSoS, two solutions were obtained and the distance to the one nearest to the correct position was recorded; the average of the results obtained in the 200 trials is also shown in Fig.~\ref{fig:mean_distance_to_true_solution} (left). As one can see, the solution of these two methods remains accurate even when the distance to the bad configuration get extremely small. The average distance between the two solutions is plotted in Fig.~\ref{fig:distance_between_slns} (left). As one can see, the 2 solutions on the sphere obtained with SoS are far apart. After the refinement of RSoS, they become nearly identical when away from the bad configuration, but remain distinct when close to it. 
Since the user is situated  directly on the sphere, the best solution is always one of the two obtained using the SoS method, as the refinement process has no sphere constraint.

In the second experiment, we proceeded similarly, except that we placed the user on a sphere of a larger radius than the earth ($0.1\%$ larger, about $6$ km) and the bad satellite position was chosen so that there are two solutions on the user sphere (rather than on earth). The SoS method was performed assuming that the user is on the earth radius and thus yielded slighly inaccurate results, as seen in Fig.~\ref{fig:mean_distance_to_true_solution}(right). RSoS gives consistently better accuracy, with two nearly identical solutions far away from the bad configuration and two distinct solutions close to it Fig.~\ref{fig:distance_between_slns}(right). ILS has a similar accuracy to RSoS far away from the bad configuration, but diverges near it.

\par{\bf Acknowledgments.} Figure~\ref{fig:Two_sats} was produced with
the help of GeoGebra. Thank you to Jan-Willem Knopper for help with
Mathematica.  

\end{document}